\documentclass[11pt]{amsart}
\usepackage{graphicx,amssymb} 
\usepackage{epsfig}

\title[Box dimension of trajectories]
{Box dimension of trajectories \\
of some discrete dynamical systems}

\author{Neven Elezovi\'c, Vesna \v Zupanovi\'c, and Darko \v Zubrini\'c }

\address{University of Zagreb, Faculty of Electrical Engineering and Computing, Unska 3, 10000 Zagreb, Croatia}

\begin{document}
\date{}

\begin{abstract}
We study the asymptotics, box dimension, and Minkowski content of trajectories of
some discrete dynamical systems. We show that under very general 
conditions,
 trajectories
corresponding to  parameters where saddle-node bifurcation appears have box dimension equal to~$1/2$,
while those corresponding to period-doubling bifurcation parameter have box dimension equal to~$2/3$.
Furthermore, all these trajectories are Minkowski nondegenerate.
The results are illustrated in the case of logistic map.
\medskip

\centerline{Contacting author: Dr. Darko \v Zubrini\'c}

\centerline{University of Zagreb, Faculty of Electrical Engineering and Computing}
\centerline{Unska 3, 10000 Zagreb, Croatia}
\medskip

\centerline{tel.\ +385 1 6129 969, fax.\ +385 1 6129 946}
\centerline{\bf darko.zubrinic@fer.hr}

\end{abstract}

\keywords
{logistic map, discrete dynamical system, box dimension, 
Minkowski content, bifurcation}

\subjclass{
37C45, 
34C23
}

\maketitle

\newtheorem{theorem}{Theorem}
\newtheorem{lemma}{Lemma}
\newtheorem{cor}{Corollary}
\newtheorem{prop}{Proposition}

\font\csc=cmcsc10

\let\ced=\c

\def\esssup{\mathop{\rm ess\,sup}}
\def\essinf{\mathop{\rm ess\,inf}}
\def\wo#1#2#3{W^{#1,#2}_0(#3)}
\def\w#1#2#3{W^{#1,#2}(#3)}
\def\wloc#1#2#3{W_{\scriptstyle loc}^{#1,#2}(#3)}
\def\osc{\mathop{\rm osc}}
\def\var{\mathop{\rm Var}}
\def\supp{\mathop{\rm supp}}
\def\Cap{{\rm Cap}}
\def\norma#1#2{\|#1\|_{#2}}

\def\C{\Gamma}

\let\text=\mbox

\catcode`\@=11
\let\ced=\c
\def\a{\alpha}
\def\b{\beta}
\def\c{\gamma}
\def\d{\delta}
\def\g{\lambda}
\def\o{\omega}
\def\q{\quad}
\def\n{\nabla}
\def\s{\sigma}
\def\div{\mathop{\rm div}}
\def\sing{{\rm Sing}\,}
\def\singg{{\rm Sing}_\ty\,}

\def\A{{\cal A}}
\def\F{{\cal F}}
\def\H{{\cal H}}
\def\W{{\bf W}}
\def\M{{\cal M}}
\def\N{{\cal N}}
\def\S{{\cal S}}

\def\eR{{\bf R}}
\def\eN{{\bf N}}
\def\Ze{{\bf Z}}
\def\Qe{{\bf Q}}
\def\Ce{{\bf C}}

\def\ty{\infty}
\def\e{\varepsilon}
\def\f{\varphi}
\def\:{{\penalty10000\hbox{\kern1mm\rm:\kern1mm}\penalty10000}}
\def\ov#1{\overline{#1}}
\def\D{\Delta}
\def\O{\Omega}
\def\pa{\partial}

\def\st{\subset}
\def\stq{\subseteq}
\def\pd#1#2{\frac{\pa#1}{\pa#2}}
\def\pdd#1#2{\frac{\pa^2#1}{\pa#2^2}}
\def\sgn{{\rm sgn}\,}
\def\sp#1#2{\langle#1,#2\rangle}

\newcount\br@j
\br@j=0
\def\q{\quad}
\def\gg #1#2{\hat G_{#1}#2(x)}
\def\inty{\int_0^{\ty}}
\def\remark{\smallskip\advance\br@j by1 \noindent{\csc Remark
\the\br@j.}\kern3mm\relax}
\def\od#1#2{\frac{d#1}{d#2}}
\def\odd#1#2{\frac{d^2#1}{d#2^2}}

\def\bg{\begin}
\def\eq{equation}
\def\bgeq{\bg{\eq}}
\def\endeq{\end{\eq}}
\def\bgeqnn{\bg{eqnarray*}}
\def\endeqnn{\end{eqnarray*}}
\def\bgeqn{\bg{eqnarray}}
\def\endeqn{\end{eqnarray}}

\def\bgeqq#1#2{\bgeqn\label{#1} #2\left\{\begin{array}{ll}}
\def\endeqq{\end{array}\right.\endeqn}

\def\abstract{\bgroup\leftskip=2\parindent\rightskip=2\parindent
        \noindent{\bf Abstract.\enspace}}
\def\endabstract{\par\egroup}

\def\udesno#1{\unskip\nobreak\hfil\penalty50\hskip1em\hbox{}
             \nobreak\hfil{#1\unskip\ignorespaces}
                 \parfillskip=\z@ \finalhyphendemerits=\z@\par
                 \parfillskip=0pt plus 1fil}
\catcode`\@=11

\def\cal{\mathcal}
\def\eR{\mathbb{R}}
\def\eN{\mathbb{N}}
\def\Ze{\mathbb{Z}}
\def\Qu{\mathbb{Q}}
\def\Ce{\mathbb{C}}

\def\sign{\mbox{\rm sign}\,}

\section{Introduction}

We are interested in bifurcation parameters
$\mu$ of discrete one-dimensional dynamical systems in the sense of nontriviality of
box dimension of the trajectory $S_\mu$, near a given trajectory of the system. 
More precisely, we are interested
in values of the parameter $\mu$ such that  $\dim_BS_\mu$ is nonzero. 
The main results are stated in Theorems~\ref{Mb} and~\ref{-Mb}.

A typical example is the
system generated by standard logistic map.
M.\ Feigenbaum studied the dynamics of
the logistic map for
$\g\in(0,4]$.
Taking $\g=\g_\ty\approx 3.570$
the corresponding invariant set $A\st[0,1]$ has both
Hausdorff and box dimensions equal to
$\approx0.538$ (Grassberger \cite{grass}, Grassberger and Procaccia \cite{grassproc}).
Here we 
compute precise values of 
box dimension of trajectories corresponding to period-doubling bifurcation parameters  $3$ and
$1+\sqrt6$, and to period-$3$ bifurcation parameter $1+\sqrt8$, see Corollary~\ref{logistic}.

Similar effect of nontriviality of box dimension of trajectories as
in bifurcation problems for discrete systems has been noticed for some planar vector fields 
having spiral trajectories $\C_\mu$, see \cite{zuzu}. There we 
have considered a standard model of Hopf-Takens bifurcation
with respect to bifurcation parameter $\mu$ where  
$\dim_B\C_\mu>1$, while  $\dim_B\C_\mu=1$ otherwise.  
We noticed that a limit cycle is born at the moment of jump of box dimension of a  spiral trajectory.
Analogously, in the case of one-dimensional discrete system a periodic trajectory is born at the moment of jump
of box dimension of a discrete trajectory (sequence).
We expect that fractal analysis of general planar spiral trajectories can be reduced to the study
of discrete one-dimensional trajectories 
 via the Poincar\'e map, see also  \cite[Remark 11]{zuzu}.
A review of results dealing with 
applications of fractal dimensions to dynamics
is given in~\cite{fdd}.

We recall the notions of box dimension and Minkowski content, see e.g.\ Mattila \cite{mattila}.
For any subset $S\st\eR^N$ by
$S_\e$ we denote the $\e$-neighbourhood of $S$ (also called Minkowski sausage of radius $\e$
around $A$, a term coined by B.\ Mandelbrot), and $|S_\e|$ is its
$N$-dimensional Lebesgue measure.
For a bounded set $S$ and given $s\ge0$ we define the upper $s$-dimensional Minkowski content of $S$ by
$$
\M^{*s}(S)=\limsup_{\e\to0}\frac{|S_\e|}{\e^{N-s}}.
$$
Analogously for the lower $s$-dimensional content of $S$. The upper box dimension of $S$ is defined by
$$
\ov\dim_BS=\inf\{s\ge0:\M^{*s}(S)=0\},
$$
and analogously the lower box dimension $\underline\dim_BS$. If both dimensions coincide, we denote it by
$\dim_BS$.
We say that a set $S$ is Minkowski nondegenerate if 
its $d$-dimensional upper and lower Minkowski contents
are in $(0,\ty)$ for some $d\ge0$, and Minkowski measurable if $\M^{*d}(S)=\M_*^d(S):=\M^d(S)\in(0,\ty)$.

Nondegeneracy of Minkowski contents for fractal strings has been
characterized by Lapidus and van Frankenhuysen \cite{lapiduspom}.
Applications of Minkowski content in the study of singular integrals can be seen in
\cite{mink} and~\cite{singl}.

For any two sequences $(a_n)_{n\ge1}$ and $(b_n)_{n\ge1}$ of positive real numbers we write $a_n\simeq b_n$ as $n\to\ty$
if there exist positive constants $A$ and $B$ such that $A\le a_n/b_n\le B$ for all $n$.
 Analogously, for two positive functions
$f,g\:(0,r)\to\eR$ we write $f(x)\simeq g(x)$ as $x\to0$ if $f(x)/g(x)\in[A,B]$ for $x$ sufficiently
small.

\section{Box dimension of some recurrently defined sequences}

The first result deals with  sequences~$(x_n)_{n\ge1}$ converging monotonically to zero.

\begin{theorem}\label{m}
Let $\a>1$ and let $f\:(0,r)\to (0,\ty)$ be a 
monotonically nondecreasing function such that $f(x)\simeq x^\a$ as $x\to0$, and
 $f(x)<x$ for all $x\in(0,r)$.
Consider the sequence $S(x_1):=(x_n)_{n\ge1}$ defined by
\bgeq
x_{n+1}=x_n-f(x_n),\q x_1\in(0,r).
\endeq
 Then
\bgeq\label{xn}
x_n\simeq n^{-1/(\a-1)}\q \mbox{as $n\to\ty$.}
\endeq
Furthermore,
\bgeq
\dim_BS(x_1)=1-\frac1\a,
\endeq
and the set $S(x_1)$ is Minkowski nondegenerate.
\end{theorem}

\begin{proof}
(a)
Assuming that $Ax^\a\le f(x)\le Bx^\a$, we have
\bgeq
0<x_{n+1}\le x_n-Ax_n^\a.
\endeq
It is easy to see that $x_n$ tends monotonically to $0$.
Using induction we first prove that $x_n\le b n^{-\b}$, where $\b:=\frac1{\a-1}$, for some positive constant $b$.
Let us consider inductive step first, and then the basis. Assume that $x_n\le b n^{-\b}$ for some $n$, and assume also
that $x_n\le x_{max}$, where $x_{max}$ is the point of maximum of $x-x^\a$, $x>0$.
Note that since $x_n$ is decreasing, converging to zero, then $x_n\le x_{max}$ for all $n$ sufficiently large.
Exploiting monotonicity of $x\mapsto x-x^\a$ on $(0,x_{max})$ we have
$$
x_{n+1}\le x_n-Ax_n^\a\le bn^{-\b}-Ab^\a n^{-\a \b}\le b(n+1)^{-\b}.
$$
In order to prove the last inequality, it suffices to show that
$$
n^{-\b}-Ab^{\a-1} n^{-\a \b}\le (n+1)^{-\b}.
$$
To this end let us consider the binomial series expansion:
\bgeqn
(n+1)^{-\b}&=&[n(1+\frac 1n)]^{-\b}=n^{-\b}+\binom{-\b}1n^{-\b-1}+R_n\label{binom}\\
&=&n^{-\b}-\b n^{-\a\b}+R_n\ge n^{-\b}-Ab^{\a-1}n^{-\a\b}.\nonumber
\endeqn
The last inequality holds provided $b$ is chosen so that $Ab^{\a-1}\ge\b$, and if $R_n\ge0$.
To prove  $R_n\ge0$ note that $R_n=a_2+a_4+a_6+\dots$, where each $a_k$ (with even $k$) has
the form
\bgeqn
a_k&=&\binom{-\b}kn^{-\b-k}+\binom{-\b}{k+1}n^{-\b-k-1}\nonumber\\
&=&\frac{\b(\b+1)\dots(\b+k-1)}{k!}n^{-\b-k}-\frac{\b(\b+1)\dots(\b+k)}{(k+1)!}n^{-\b-k-1}.\nonumber
\endeqn
Inequality $a_k\ge0$ is equivalent with
$$
n\ge\frac{\b+k}{k+1}=\frac1{(k+1)(\a-1)}+\frac{k}{k+1}.
$$
For all even $k$ the right-hand side obviously does not exceed
$n_0=n_0(\a)=\frac1{3(\a-1)}+1$.
The condition $x_n\le x_{max}$ for $n\ge n_0$ is secured if we take $n_0$ sufficiently large.
From condition $Ab^{\a-1}\ge\b$ we see that we must take $b\ge (\b/A)^{\b}$. Hence, the basis of induction and
inductive step hold for $n\ge n_0$ with such a $b$. Taking $b$ still larger, we 
can achieve that $x_n\le bn^{-\b}$ for all $n\ge1$.

(b) To prove that there exists $a>0$ such that $x_n\ge an^{-\b}$ for all $n\ge1$, 
we use only $x_{n+1}\ge x_n-Bx_n^\a$. 
Assuming by induction that the desired inequality holds for a fixed $n$ we obtain analogously as in
(a) that
\bgeq\label{ineq}
x_{n+1}\ge x_n-Bx_n^\a\ge an^{-\b}-Ba^\a n^{-\a \b}\ge a(n+1)^{-\b},
\endeq
under the assumption that $x_n\le x_{max}$.
In order to show the last inequality in (\ref{ineq})
we use binomial expansion  (\ref{binom}) again, and proceed by writing $\b=\c+\d$ with arbitrarily chosen
positive constants $\c$ and $\d$. We have to achieve
\bgeqn
(n+1)^{-\b}&=&(n^{-\b}-\c n^{-\a\b})-(\d n^{-\a\b}-R_n)\nonumber\\
&\le& n^{-\b}-B a^{\a-1}n^{-\a\b}.\nonumber
\endeqn
This holds provided $\c\ge Ba^{\a-1}$, that is, $a\le(\c/B)^\b$, and if $R_n\le \d n^{-\a\b}$ for some
$\d\in(0,\b)$. Note that $R_n\le \d n^{-\a\b}$  is equivalent with
$$
n\sum_{k=2}^\ty\binom{\b}kn^{-k}\le \d.
$$
that is, with
$$
n\left[(1+\frac1n)^\b-1-\binom\b 1\frac1n\right]\le\d.
$$
Using Taylor's formula $(1+\frac1n)^\b=1+\binom \b1\frac1n+\binom\b2\ov x^2$, $0<\ov x<\frac1n$,
we see that the above inequality is satisfied when $\binom\b2n\ov x^2\le\d$, that is, when
$\binom\b2n^{-1}\le\d$.
This holds for all $n\ge n_0$ if $n_0$  is large enough.
We can choose $n_0$ large enough so that also $x_{n_0}\le x_{max}$.
Taking $a$ small enough we can achieve the basis of induction, $x_{n_0}\ge a n_0^{-\b}$.
Taking $a>0$ still smaller 
the lower bound will hold for all $n\ge1$. 
This completes the proof of the lower bound of $x_n$ by induction.

(c) Since $f$ is nondecreasing, the sequence $l_n:=x_n-x_{n+1}=f(x_n)$ is nonincreasing.
Hence, we can derive Minkowski nondegeneracy of $S(x_1)$ using
Lapidus and Pomerance \cite[Theorem 2.4]{lapiduspom}.
Indeed,  from (\ref{xn}) we have
$$
l_n=f(x_n)\simeq x_n^\a\simeq n^{-\a/(\a-1)}=n^{-1/d},
$$
where $d:=1-\frac1\a\in(0,1)$.
Using the mentioned result 
we conclude that $S(x_1)$ is Minkowski nondegenerate and $\dim_BS(x_1)=d=1-\frac1\a$.
\end{proof}

\remark Step (c) in the proof of Theorem~\ref{m}
can be carried out by directly estimating Minkowski contents of $S=S(x_1)$.
Using $l_n:=x_n-x_{n+1}=f(x_n)\le B(bn^{-\b})^\a=Bb^\a n^{-\a\b}$ we see that
$l_n\le2\e$ if $n\ge(\frac12 Bb^\a)^{d}\e^{-d}$, where $d:=1-\frac1\a$.
Defining $n_0=n_0(\e):=\lceil
(\frac12 Bb^\a)^{d}\e^{-d}
\rceil$
we have
\bgeqn\label{ge}
|S_\e|\ge x_{n_0}+2\e(n_1-1),
\endeqn
where $n_1=n_1(\e)$ is obtained in the similar way from the condition 
$l_n=f(x_n)\ge2\e$. It is satisfied for $n\le n_1:=\lfloor (\frac12Aa^{\a})^d\e^{-d}\rfloor$.
Using (\ref{ge}) we conclude that
\bgeq\label{gm}
\M_*^d(S)\ge\frac ab\left(\frac 2{B}\right)^{1/\a}+2\left(\frac12Aa^\a\right)^d.
\endeq
In the analogous way, estimating $|S_\e|$ from above, we obtain
\bgeq\label{dm}
\M^{*d}(S)\le\frac ba\left(\frac 2{A}\right)^{1/\a}+2\left(\frac12Bb^\a\right)^d.
\endeq
This proves that $S$ is Minkowski nondegenerate and $\dim_BS=d$.

\remark  We do not know if the set $S=S(x_1)$ corresponding to $f(x)= A\cdot x^\a$ in Theorem \ref{m},
where $A>0$, is Minkowski measurable.
Numerical experiments show that in this case any corresponding sequence $S=(x_n)_{n\ge1}$, $x_1\in(0,1)$, 
is Minkowski measurable, and
\bgeq
\M^d(S)=\left(\frac2{A}\right)^{1/\a}\frac\a{\a-1}.
\endeq
This value is obtained if we let formally $a=b=(\b/A)^\b$ in (\ref{gm}) and (\ref{dm}).

\medskip

The following result deals with sequences $(x_n)_{n\ge1}$ oscillating around a fixed point $x_0$, so that
their two subsequences defined by odd and even indices monotonically converge to~$x_0$.
It suffices to consider the case $x_0=0$.

\begin{theorem}\label{-m}
Let $f:(-r,r)\to\eR$ be a function
such that $f(x)\simeq |x|^\a$ as $x\to0$, where $\a>1$.
We also assume that the function $f(-x-f(x))-f(x)$ satisfies the following conditions:
\bgeqn
&\mbox{it is monotonically nondecreasing for $x>0$ small enough,}&\label{mon}\\
&\mbox{it is monotonically nonincreasing for $x<0$ small enough,}&\nonumber\\
&f(-x-f(x))-f(x)\simeq\pm|x|^{2\a-1}\q\mbox{as $x\to0\pm$.}\label{funkc}&
\endeqn
Then there exists $r_1>0$ such that for any sequence $S(x_1):=(x_n)_{n\ge1}$ defined by
\bgeq
x_{n+1}=-x_n-f(x_n),\q x_1\in(-r_1,r_1),
\endeq
we have
\bgeq\label{xnn}
|x_n|\simeq n^{-1/(2\a-2)},\q \mbox{as $n\to\ty$.}
\endeq
Furthermore,
\bgeq
\dim_BS(x_1)=1-\frac1{2\a-1},
\endeq
and the set $S(x_1)$ is Minkowski nondegenerate.
\end{theorem}

\begin{proof}
Note that if we define $F(x):=-x-f(x)$ then $F^2(x)=F(F(x))=x-g(x)$, where $g(x):=f(-x-f(x))-f(x)$. 
We have
\bgeq\label{gx}
g(x)\simeq \pm|x|^{2\a-1}\q\mbox{as $x\to0\pm$.}
\endeq
As $2\a-1>1$ and $g(0)=0$, from (\ref{gx}) we see that there exists $r_1\in(0,r)$ such that
$0<g(x)<x$ for $x\in(0,r_1)$ and $-x<g(x)<0$ for $x\in(-r_1,0)$.

Starting with  $x_1\in(0,r_1)$, the sequence $y_n:=x_{2n-1}$, $n\ge1$, 
satisfies $y_{n+1}=y_{n}-g(y_{n})$, and since $(y_n)$ is nonincreasing, it is contained in $(0,r_1)$.
By Theorem \ref{m} applied to $g$ and the sequence $(y_n)$
we have that $y_n=x_{2n-1}\simeq n^{-1/(2\a-2)}$ as $n\to\ty$. To obtain the same asymptotics for $|x_{2n}|$, it suffices
to start with $x_2=F(x_1)<0$, and to consider the sequence $z_n:=x_{2n}$,
${n\ge1}$, contained in $(-r,0)$. Using again Theorem 
\ref{m} (modified to this situation; note that $-x<g(x)<0$) 
applied to the sequence $z_n$, we obtain $|z_n|=|x_{2n}|\simeq n^{-1/(2\a-2)}$. Hence $|x_{n}|\simeq n^{-1/
(2\a-2)}$. The same asymptotics is obtained if we start with $x_1\in(-r_1,0)$.

Exploiting finite stability of the upper box dimension, see Falconer \cite[p.\ 44]{falc},
we have that $\ov\dim_{B}S=\max\{\ov\dim_BS_-,\ov\dim_BS_+\}$, where $S_-:=S\cap(-r,0)$ and
$S_+:=S\cap(0,r)$ are negative and positive part of the sequence $S=S(x_1)$. 
Since $\ov\dim_BS_-=\ov\dim_BS_+=1-\frac1{2\a-1}$, see Theorem \ref{m}, we conclude that
$\ov\dim_BS=1-\frac1{2\a-1}$.

To estimate the lower box dimension, first note that the sets 
$S_{+}$ and $S_-$
are separated by $x=0$, hence $(S_+)_\e\cap (S_-)_\e=(-\e,\e)$. 
Therefore $|S_\e|=|(S_+)_\e|+|(S_-)_\e|-|(S_+)_\e\cap(S_-)_\e|=|(S_+)_\e|+|(S_-)_\e|-2\e$, and
from this we immediately obtain that
$\M_*^d(S)\ge\M_*^d(S_+)+\M_*^d(S_-)>0$, where $d:=1-\frac1{2\a-1}$.
We have used also Minkowski nondegeneracy of $S_{\pm}$. Hence,
$\underline\dim_BS\ge d$. This finishes the proof of $\dim_BS=d$.
\end{proof}

\remark
Conditions of Theorem \ref{-m} are satisfied when for example $f(x)=|x|^\a$, $x\in(-r,r)$.

In order to facilitate the study of bifurcation problems below,
it will be convenient to formulate the following consequences of Theorems~\ref{m} and~\ref{-m}.

\begin{theorem}\label{M} 
Let $F:(x_0-r,x_0+r)\to\eR$ be a function of class $C^3$, such that
\bgeqn
F(x_0)&=&x_0,\label{0}\\
 F'(x_0)&=&1,\label{od}\\
F''(x_0)&<&0.\label{odd}
\endeqn
Then there exists $r_1>0$ such that for any sequence $S(x_1)=(x_n)_{n\ge1}$ defined by
$$
x_{n+1}=F(x_n),\q x_1\in(x_0,x_0+r_1),
$$
we have $|x_n-x_0|\simeq n^{-1}$ as $n\to\ty$,
\bgeq
\dim_BS(x_1)=\frac12,
\endeq
and $S(x_1)$ is Minkowski nondegenerate.
Analogous result holds if $x_1\in(x_0-r_1,x_0)$, assuming that in (\ref{odd})
we have the opposite sign.
\end{theorem}

\begin{proof}
We can assume without loss of generality that $x_0=0$, and let $x>0$. By the Taylor formula,
using (\ref{0}), (\ref{od}), we have that
$$
F(x)=x-f(x),
$$
where
$$
f(x)=-\frac{F''(0)}{2!}x^2-\frac{F'''(\ov x)}{3!}x^3
$$
with $\ov x=\ov x(x)\in(0,r)$. Since $F''(0)<0$, we see that $f(x)\simeq x^2$. The condition $f(x)<x$
is clearly satisfied in $(0,r_1)$ for $r_1$ sufficiently small. The function $f$ is increasing,
since  using again Taylor's formula applied to $F'\in C^2$, we get
$$
f'(x)=1-F'(x)=-F''(0)\,x-\frac{F'''(\ov x)}{2!}\,x^2>0
$$
for $x\in(0,r_1)$ with $r_1$ small enough.
The claim follows from Theorem~\ref{m} with $\a=2$.
Analogously for $x\in(-r_1,0)$.
\end{proof}

\begin{theorem}\label{-M}
Let $F:(x_0-r,x_0+r)\to\eR$ be a function of class $C^4$, such that
\bgeqn
F(x_0)&=&x_0,\label{0-}\\
 F'(x_0)&=&-1,\label{-od}\\
 F''(x_0)&\ne&0.\label{-odd}
\endeqn
Then there exists $r_1>0$ such that for any sequence $S(x_1)=(x_n)_{n\ge1}$ defined by
$$
x_{n+1}=F(x_n),\q x_1\in(x_0-r_1,x_0+r_1),\q x_1\ne x_0,
$$
we have $|x_n-x_0|\simeq n^{-1/2}$ as $n\to\ty$,
\bgeq
\dim_BS(x_1)=\frac 23,
\endeq
and $S(x_1)$ is Minkowski nondegenerate.
\end{theorem}

\begin{proof}
We assume without loss of generality that $x_0=0$. It suffices to check that conditions
of Theorem~\ref{-m} are satisfied. Using Taylor's formula, (\ref{0}), and (\ref{-odd}), we get
$$
F(x)=-x-f(x),
$$
where
$$
f(x)=-\frac{F''(0)}{2!}\,x^2-\frac{F^{'''}(\ov x)}{3!}\,x^3,
$$
where $\ov x=\ov x(x)\in(-r,r)$. 
Now we consider the function $g(x)=f(-x-f(x))-f(x)$:
$$
g(x)=-F''(0)x\cdot f(x)+\dots=\frac12F''(0)^2x^3+\mbox{higher order terms}.
$$
Since $g'(x)=\frac32F''(0)^2\,x^2+\mbox{higher order terms}>0$ for $x\ne0$ such that $|x|$ is small enough,
it is clear that $g(x)$ is increasing on $(-r_1,r_1)$, provided $r_1>0$ is small enough.
Also,
$$
g(x)=\frac{F''(0)^2}{2}x^3+o(x^3)\simeq x^3\q\mbox{as $x\to0$.}
$$
This shows that conditions (\ref{mon}) and (\ref{funkc}) are fulfilled.
The claim follows from Theorem~\ref{-m} with $\a=2$.
\end{proof}

\remark It is clear that more general versions of Theorems~\ref{M} and~\ref{-M} can be proved where
finitely many consecutive derivatives of $F$ of orders $k=2,3,\dots,2m-1$ are equal to zero,
and $F^{(2m)}(x_0)\ne0$.

\begin{lemma}\label{hyperbolic} 
Assume that $S=(x_n)\st\eR$ is a sequence of positive numbers converging exponentially to zero, that is,
there exists $\g\in(0,1)$ and a constant $C>0$ such that $0<x_n\le C\g^n$ for all $n$. Then
$\dim_BS=0$.
\end{lemma}

\begin{proof} For any fixed $\e>0$ inequality $x_n\le \g^n<2\e$ is satisfied for $n\ge n_0(\e):=
\lceil\frac{\log(2\e)}{\log\g}\rceil$. Hence, 
$$
|S_\e|\le2\e+n_0(\e)\cdot 2\e,
$$
and from this $\M^{*s}(S)=0$ for any $s\in(0,1)$. Therefore, $\ov\dim_B S=0$.
\end{proof}

From this we immediately obtain the following result. 
The notion of hyperbolic fixed point of the system is described e.g.\ in
Devaney~\cite{devaney}.

\begin{theorem} {\rm (Hyperbolic fixed point)}\label{hyp}
Let $F:(x_0-r,x_0+r)\to\eR$ be a function of class $C^1$,
$F(x_0)=x_0$, and $|F'(x_0)|<1$. Then there exists $r_1\in(0,r)$
such that for any sequence $S(x_1)=(x_n)_{n\ge1}$ defined by
$$
x_{n+1}=F(x_n),\q x_1\in(x_0-r_1,x_0+r_1)
$$
we have
$$
\dim_BS(x_1)=0.
$$
\end{theorem}

It is easy to see that under the assumptions on $f$ given in Theorem
\ref{m} we have that for each $\g\in(0,1)$ the sequence $S:=(x_n)_{n\ge1}$ 
corresponding to $x_{n+1}=\g x_n-g(x_n)$, $x_1\in(0,r_\g)$, with $r_\g$ sufficiently small, has exponential decay,
$0<x_n\le \g^n$.
Hence
$\dim_BS=0$.
\smallskip

Now we state a simple but useful comparison result, which we shall need in the proof of Theorem~\ref{exp}.

\begin{lemma}\label{comp} {\rm (Comparison principle for box dimensions)}
Assume that $A=(a_n)_{n\ge1}$ and $B=(b_n)_{n\ge1}$ are two decreasing sequences of positive real numbers
converging to zero, such that
the sequences of their differences $(a_n-a_{n+1})_{n\ge1}$
and $(b_n-b_{n+1})_{n\ge1}$ are monotonically nonincreasing. If $a_n\le b_n$ then
$$
\ov\dim_BA\le\ov\dim_BB,\q
\underline\dim_BA\le\underline\dim_BB.
$$
\end{lemma}

\begin{proof}
Using the fact that the Borel rarefaction index of $A$ is equal to the upper box dimension, 
see Tricot \cite[p.\ 34 and Theorem on p.\ 35]{tricot}, we obtain
\bgeq
\ov\dim_BA=\underset{n\to\ty}{\ov\lim}\,\,\frac1{1+\frac{\log(1/a_n)}{\log n}}.
\endeq
Since $0<a_n\le b_n$ we conclude that $\ov\dim_B A\le\ov\dim_B B$.
Using methods described in Tricot \cite[pp.\ 33--36]{tricot} it can be shown that 
analogous result holds for the lower box dimension:
\bgeq\label{diml}
\underline\dim_BA=\underset{n\to\ty}{\underline\lim}\,\,\frac1{1+\frac{\log(1/a_n)}{\log n}}.
\endeq
This immediately implies $\underline\dim_B A\le\underline\dim_B B$.
\end{proof}

Now we consider a sequence with a very slow convergence to $0$, such that its box dimension is maximal possible.
We achieve this by assuming that $f(x)$ converges very fast to $0$ as $x\to0$. An example of such a function
is $f(x)=\exp(-1/x)$.

\begin{theorem}\label{exp}
Let $f\:(0,r)\to (0,\ty)$ be a nondecreasing function such that $f(x)<x$ and for any $\a>1$ we have
 that $f(x)=O(x^\a)$ as  $x\to0$.
Consider the sequence $S:=(x_n)_{n\ge1}$ defined by
$x_{n+1}=x_n-f(x_n)$, $x_1\in(0,r)$.
Then
\bgeq
\dim_BS=1.
\endeq
\end{theorem}

\begin{proof} It is clear that $x_n\to0$. 
For any fixed $\a>1$ there exists $B_\a>0$ such that we have $f(x)\le B_\a x^\a$,
hence $x_{n+1}\ge x_n-B_\a x_n^\a$ for all $n\ge1$. As in step (b) 
of the proof of Theorem~\ref{m}
we conclude that there exists $a=a_\a>0$
such that $x_n\ge an^{-1/(\a-1)}$ for all $n$. 
Since $x_n\to0$ monotonically, then $c_n=x_n-x_{n+1}=f(x_n)\to0$ also monotonically.
Therefore, using Lemma~\ref{comp} (see also (\ref{diml})), we get 
$$
\underline\dim_BS\ge\underline\dim_B\{an^{-1/(\a-1)}\}=\frac1{1+(\a-1)^{-1}}= 1-\frac1\a.
$$ 
The claim follows
by letting $\a\to\ty$.
\end{proof}

\section{Box dimension of trajectories at bifurcation points}

For definitions of saddle-node (or tangent) bifurcation and period-doubling bifurcation
and basic results see Devaney \cite[Section 1.12]{devaney}.
Note that conditions in Theorem~\ref{Mb} are essentially the same as those in \cite[Theorem~12.6]{devaney}.
Also conditions in Theorem~\ref{-Mb} are essentially the same as those in \cite[Theorem~12.7]{devaney}.
The novelty in Theorems~\ref{Mb} and~\ref{-Mb} are precise values of box dimensions of trajectories near the point
of bifurcation, convergence rate of trajectories, and their Minkowski nondegeneracy.

\begin{theorem}\label{Mb} {\rm (Saddle-node bifurcation)}
Suppose that a function $F:J\times(x_0-r,x_0+r)\to\eR$, where $J$ is an open interval in $\eR$,
is such that $F(\g_0,\cdot)$ is of class $C^3$ for some $\g_0\in\eR$, 
and $F(\cdot,x)$ of class $C^1$ for all~$x$. Assume that
\bgeqn
F(\g_0,x_0)&=&x_0,\\
\pd F{x}(\g_0,x_0)&=&1,\label{odb}\\
\pdd Fx(\g_0,x_0)&<&0,\label{oddb}\\
\pd F{\g}(\g_0,x_0)&\ne&0.
\endeqn
Then $\g_0$  is the point where saddle-node  bifurcation occurs. Furthermore, there exists $r_1\in(0,r)$ such that 
for any sequence $S(\g_0,x_1)=(x_n)_{n\ge1}$ defined by
$$
x_{n+1}=F(\g_0,x_n),\q x_1\in(x_0,x_0+r_1),
$$
we have $|x_n-x_0|\simeq n^{-1}$ as $n\to\ty$,
\bgeq
\dim_BS(\g_0,x_1)=\frac12,
\endeq
and $S(\g_0,x_1)$ is Minkowski nondegenerate.
Analogous result holds if $x_1\in(x_0-r_1,x_0)$, assuming that in (\ref{oddb})
we have the opposite sign.
\end{theorem}

\begin{proof}
The claim follows from Theorem~\ref{M} and \cite[Theorem~12.6]{devaney}.
\end{proof}

\begin{theorem}\label{-Mb} {\rm(Period-doubling bifurcation)}
Let $F:J\times(x_0-r,x_0+r)\to\eR$ be a function of class $C^2$, 
where $J$ is an open interval in $\eR$,
and $F(\g_0,\cdot)$ is of class $C^4$ for some $\g_0\in J$.
 Assume that
\bgeqn
F(\g_0,x_0)&=&x_0,\\
\pd Fx(\g_0,x_0)&=&-1,\label{-odb}\\
\pdd Fx(\g_0,x_0)&\ne&0,\label{-oddb}\\
\frac{\partial^2(F^2)}{\partial\g\,\partial x}(\g_0,x_0)\ne0,&&\frac{\partial^3(F^2)}{\partial x^3}
(\g_0,x_0)\ne0,
\endeqn
where we have denoted $F^2=F\circ F$. Then $\g_0$ is the point where
period-doubling bifurcation occurs. Furthermore,
there exists $r_1\in(0,r)$ such that for any sequence $S(\g_0,x_1)=(x_n)_{n\ge1}$ defined by
$$
x_{n+1}=F(\g_0,x_n),\q x_1\in(x_0-r_1,x_0+r_1),\q x_1\ne x_0,
$$
we have $|x_n-x_0|\simeq n^{-1/2}$ as $n\to\ty$,
\bgeq
\dim_BS(\g_0,x_1)=\frac 23,
\endeq
and $S(\g_0,x_1)$ is Minkowski nondegenerate.
\end{theorem}

\begin{proof}
The claim follows from Theorem~\ref{-M} and \cite[Theorem~12.7]{devaney}.
\end{proof}

Now we apply preceding results to bifurcation problem generated by  the logistic map.
By $d(x,A)$, where $x\in\eR$ and $A\st\eR$,  we denote Euclidean distance from $x$ to $A$.

\begin{cor}\label{logistic} {\rm(Logistic map)} Let $F(\g,x)=\g x(1-x)$, $x\in(0,1)$, and let
 $S(\g,x_1)=(x_n)_{n\ge1}$ be a sequence defined by initial value $x_1$ and $x_{n+1}=F(\g,x_n)$.

(a) For $\g_0=1$, taking $x_1>0$ sufficiently close to $x_0=0$, we have that
$x_n\simeq n^{-1}$ as $n\to\ty$, and
$$
\dim_BS(1,x_1)=\frac12.
$$

(b) (Onset of period-$2$ cycle) For $\g_0=3$
the corresponding fixed point is $x_0=2/3$.
For any $x_1$ sufficiently close to $x_0$ we have that $|x_n-x_0|\simeq n^{-1/2}$, and
$$
\dim_BS(3,x_1)=\frac23.
$$

(c) For any $\g\notin\{1,3\}$ and $x_1$ such that the sequence $S(\g,x_1)$ is convergent,
we have that $\dim_BS(\g,x_1)=0$.

(d) (Onset of period-$4$ cycle)  If $\g_0=1+\sqrt6$ then
for any $x_1$ sufficiently close to period-$2$  trajectory $A=\{a_1,a_2\}$ we have that
 $d(x_n,A)\simeq n^{-1/2}$ as $n\to\ty$, and
$$
\dim_BS(1+\sqrt6,x_1)=\frac23.
$$

(e) (Period-$3$ cycle) Let $\g_0=1+\sqrt8$ and let $a_1$, $a_2$, $a_3$ be fixed points of  $F^3$
such that $0<a_1<a_2<a_3<1$, $F(a_1)=a_2$, $F(a_2)=a_3$, and $F(a_3)=a_1$. 
Then there exists $\d>0$ such that for any initial value
$$
x_1\in(a_1-\d,a_1)\cup(a_2-\d,a_2)\cup(a_3,a_3+\d)
$$
we have $d(x_n,\{a_1,a_2,a_3\})\simeq n^{-1}$
as $n\to\ty$,
and
$$
\dim_BS(1+\sqrt8,x_1)=\frac12.
$$
All trajectories appearing in this corollary are Minkowski nondegenerate.
\end{cor}

\begin{proof} Claim (a) follows from Theorem~\ref{m}. For (b) and (c)  see Theorems~\ref{-Mb} and \ref{hyp}.
Claim (d) follows from Theorem~\ref{-M}  since $(F^2)''(\g_0,x_0)\ne0$.  

Claim (e) follows using Theorem~\ref{M} applied to $F^3$. Note that these three intervals are
disjoint for $\d>0$ small enough. See Strogatz \cite[pp.\ 362 and 363]{strogatz}.
The fact that for $\g_0=1+\sqrt8$ we have $(F^3)''(a_i)\ne0$, $i=1,2,3$,  can be obtained 
by direct computation.
\end{proof}


\medskip

\remark It would be interesting to know precise values of box dimensions of trajectories corresponding
to all period-doubling bifurcation parameters $\g_k$ where $2^k$-periodic points occur.
On the basis of numerical experiments we expect that all of them will be equal to $2/3$.


\medskip

{\csc Example.}  For $F(\g,x):=\g e^{x}$, see Devaney \cite[Section 1.12]{devaney}, we can obtain similar results.
Indeed, using Theorem~\ref{Mb}  we obtain that
$$
\dim_BS(e^{-1},x_1)=\frac12
$$ 
for all $x_1$ in a punctured neighbourhood of $x_0=1$.
Using Theorem~\ref{-Mb} we obtain that for any $x_1$ in a punctured neighbourhood of $x_0=-1$ we have
$$
\dim_B(-e,x_1)=\frac23.
$$
If $\g\notin\{e^{-1},-e\}$, we have $\dim_BS(\g,x_1)=0$ provided
$S(\g,x_1)$ is a convergent sequence, see Theorem~\ref{hyp}.


\begin{thebibliography}{10}

\let\small=\rm
\let\normalsize=\rm


\bibitem{devaney}  Devaney R.L., {\em An Introduction to Chaotic Dynamical Systems}, The Benjamin/Cummings, New York, 1986.

\bibitem{falc} Falconer K., {\em Fractal Geometry}, John Wiley and Sons, Chichester
1990.


\bibitem{grass}
Grassberger, P., On the Hausdorff dimension of fractal attractors, J.\ Statist.\ Phys., 26 (1981), 1, 173--179.

\bibitem{grassproc} 
Grassberger P., Procaccia I., Measuring the strangeness of strange attractors, Phys.\ D 9 (1983),  1-2, 189--208.







\bibitem{lapiduspom} Lapidus M.L., Pomerance C., 
The Riemann zeta-function and the one-dimensional Weyl-Berry conjecture for fractal drums,
Proc.\ London Math.\ Soc.\ (3) {\bf 66} (1993), no. 1, 41--69.

\bibitem{mattila} Mattila P., {\it Geometry of Sets and Measures in Euclidean Spaces. Fractals and
Rectifiability}, Cambridge 1995.

\bibitem{strogatz} Strogatz S.H., 
{\it Nonlinear Dynamics and Chaos: With Applications to Physics, Biology, Chemistry and Engineering},
Addison Wesley, 1994.



\bibitem{tricot} Tricot C., {\em Curves and Fractal Dimension}, Springer--Verlag, 1995.

\bibitem{mink}
\v Zubrini\'c D.,
Minkowski content and singular integrals, {\it Chaos,
Solitons and Fractals}, {\bf 17/1} (2003), 169--177.

\bibitem{singl} \v Zubrini\'c D., Singular sets of Lebesgue integrable functions,
{\it Chaos, Solitons and Fractals}, {\bf 21} (2004) 1281--1287.


\bibitem{zuzu} \v Zubrini\'c D., \v Zupanovi\'c V., Fractal analysis of spiral trajectories of some
planar vector fields, Bulletin des Sciences Math\'ematiques, 129/6 (2005), 457-485.


\bibitem{fdd} \v Zupanovi\'c V., \v Zubrini\'c D., Fractal dimensions in dynamics, 
 in {\it Encyclopedia of Mathematical Physics},  
Jean-Pierre Fran\ced{c}oise, Greg Naber, Sheung Tsun Tsou (editors), Elsevier Academic Press, 2006, to appear.


\end{thebibliography}
\end{document}